\documentclass[12pt,reqno]{amsart}
\usepackage{amssymb}
\usepackage[usenames,dvipsnames]{color}
\usepackage{euscript}
\usepackage{multicol}
\usepackage{mathrsfs}
\usepackage{amssymb}\usepackage{amsmath}
\usepackage{times} 
\usepackage{graphicx}

\newcommand{\be}[1]{\begin{equation}\label{#1}}
\newcommand{\ee}{\end{equation}}
\newtheorem{theorem}{Theorem}

\newtheorem{remark}{Remark}

\newtheorem{lemma}{Lemma}
\newtheorem{proposition}{Proposition}

\numberwithin{equation}{section}

\begin{document}

\title[Hardy and Sobolev inequalities on antisymmetric functions]
{Hardy and Sobolev inequalities on antisymmetric functions}

\author{Th. Hoffmann-Ostenhof}
\address{Thomas Hoffmann-Ostenhof: University of Vienna}
\email{thoffmann@tbi.univie.ac.at}

\author{A. Laptev}
\address{Ari Laptev: Imperial College London \\ 180 Queen's Gate \\ London SW7 2AZ \\ UK  
and Sirius Mathematics Center, Sirius University of Science and Technology,
1 Olympic Ave, 354340, Sochi, Russia}
\email{a.laptev@imperial.ac.uk}

\author{I. Shcherbakov}
\address{Ilya Shcherbakov: Sirius Mathematics Center
Sirius University of Science and Technology,
1 Olympic Ave, 354340, Sochi, Russia} 
\email{stscherbakov99@yandex.ru}  



\begin{abstract}
We obtain sharp Hardy inequalities on  antisymmetric functions where antisymmetry is understood for multi-dimensional particles. Partially it is an extension of the previously published paper \cite{HL}, where Hardy's inequalities were considered for the antisymmetric functions  in the case of the 1D particles.
As a byproduct we obtain some Sobolev and  Gagliardo-Nirenberg type  inequalities that are applied to the study of spectral properties of Schr\"odinger operators.
\end{abstract}

\date{}

\maketitle

\section{Introduction}

\noindent
The classical Hardy inequality reads for $u\in \mathcal H^1(\Bbb R^n)$, $n\ge 3$,
\begin{equation}\label{Hclass}
 \int_{\Bbb R^n}|\nabla u(x)|^2dx\ge \frac{(n-2)^2}{4}\int_{\Bbb R^n}\frac{|u(x)|^2}{|x|^2}dx.
\end{equation}
The literature concerning different versions of Hardy inequalities and their applications is extensive and we are not able to cover it in this paper. We just mention the classical paper \cite{B} and books \cite{BEL1},  \cite{D1}, \cite{D2},   \cite{FLW}, \cite{M}. 
Clearly if $n=2$, then the Hardy inequality \eqref{Hclass} is trivial. 

\medskip
\noindent
For $n\ge 3$ we also have the classical Sobolev inequality 
\begin{equation}\label{SobClass}
\int_{\mathbb{R}^n} |\nabla u(x) |^2 dx \ge S(n) \left(\int_{\mathbb{R}^n}|u(x)|^{\frac{2n}{n - 2}} \, dx\right)^{\frac{n-2}{n}}, \quad u\in \mathcal H^1(\Bbb R^n),
\end{equation}
where 
\begin{equation}\label{SobConst}
S(n) = \frac{n(n-2)}{4} \, |\Bbb S^{n}|^{2/n} =  \frac{n(n-2)}{4} \, 2^{2/n} \, \pi^{1+1/n} \, 
\Gamma\left(\frac{n+1}{2}\right)^{-2/n}.
\end{equation}
The inequalities  \eqref{Hclass} and \eqref{SobClass} are related. It is known, see \cite{BEL1}  and \cite{Sei},  that \eqref{Hclass} implies \eqref{SobClass} and the Sobolev inequality \eqref{SobClass} implies a weak version of the Hardy's inequality.
Hardy and Sobolev inequalities are also closely related to spectral properties of the negative eigenvalues of Schr\"odinger operators. 

\medskip
\noindent
The aim of this paper is to consider functional and spectral inequalities and their relations on antisymmetric functions. 

\medskip
\noindent
Let $N$ and $d$ be natural numbers. We consider $x = (x_1, \dots, x_N) \in \mathbb{R}^{dN},$ where 
$x_i = (x_{i1},\dots, x_{id}) \in \mathbb{R}^d$ for all $1\le i \le N.$ Every function $u$ defined on $\mathbb{R}^{dN}$ we call \emph{antisymmetric} hereafter, if for all $1 \le i,j \le N$ and $x_1, \dots, x_N \in \mathbb{R}^d$
$$
u(x_1,\dots,x_i,\dots,x_j,\dots,x_N) = - u(x_1,\dots,x_j,\dots,x_i,\dots,x_N).
$$
Let us consider the subclass  of antisymmetric functions  from $\mathcal H^1(\Bbb R^{dN})$, that we denote by $\mathcal H_A^1(\Bbb R^{dN})$. 
 Clearly, $\mathcal H^1_A(\Bbb R^{dN}) \subset \mathcal H^1(\Bbb R^{dN})$ and therefore it is expected that the constants in the inequalities \eqref{Hclass} and \eqref{SobClass} are larger.

\medskip
\noindent 
Let $\mathcal{V}_d(N)$ be the degree of the Vandermonde determinant defined in \eqref{det}.
Among the main results obtained in this paper is the following theorem.
\begin{theorem}\label{Hardy_A}
 Let $d,N\in \Bbb N$ and let $d\ge1$, $N\ge2$. For any $u\in \mathcal H^1_A (\Bbb R^{dN})$ we have 
\begin{equation}\label{H_A}
 \int_{\Bbb R^{dN} }|\nabla u(x)|^2 \, dx\ge H_A(dN)\, \int_{\Bbb R^{dN}}\frac{|u(x)|^2}{|x|^2}\, dx,
\end{equation}
where 
\begin{equation}\label{hardydn}
H_A(dN) = \frac {(dN-2)^2}4 + \mathcal{V}_d(N)(\mathcal{V}_d(N) + dN - 2).
\end{equation}
\end{theorem}
We shall show at the end of Section \ref{HardyIn}  that the constant \eqref{hardydn}  in \eqref{H_A} is sharp. The constant $(dN-2)^2/4$ in \eqref{hardydn} is the classical Hardy constant and the constant $\mathcal{V}_d(N)(\mathcal{V}_d(N) + dN - 2)$ is related to the antisymmetry of the class of functions $\mathcal H^1_A (\Bbb R^{dN})$.

\medskip
\noindent
The next theorem gives an improvement of the constant in the Sobolev inequality restricted to antisymmetric functions.

\begin{theorem}\label{Sobolev_A}
For every $u\in \mathcal H_A^1(\mathbb{R}^{dN}),$ $dN \ge 3$,
\[
\int_{\mathbb{R}^{dN}}|\nabla u(x)|^2 dx \ge S_A(dN) \left(\int_{\mathbb{R}^{dN}}|u(x)|^{\frac{2dN}{dN-2}}dx\right)^{\frac{dN-2}{dN}},
\]
where 
$$
S_A(dN) = (N!)^{\frac 2{dN}} S(dN) = 
\pi dN(dN-2)\left(\frac{\Gamma(\frac{dN}2)N!}{ \Gamma(dN)}\right)^{\frac 2{dN}}
$$ 
and where $S(dN)$ is the classical Sobolev constant in $\Bbb R^{dN}$.
\end{theorem}

The constant $S_A(dN)$ is sharp and substantially larger than the classical constant $S(dN)$.

 \medskip
 \noindent
 Note that the inequality  \eqref{H_A} in the case $d=1$ and arbitrary $N$ has been obtained in \cite{HL} with the sharp constant
 $H_A(N) = \frac{(N^2-2)^2}{4}$, $N\ge2$. Comparing this with the classical constant $(N-2)^2/4 \sim N^2$ as $N\to\infty$, the constant  $H_A(N)\sim N^4$. The proof of Theorem \ref{Hardy_A} in the case $d=1$ is based on the lowest eigenvalue of the Laplace-Beltrami operator on spherical harmonics generated by antisymmetric harmonic polynomials.
The case $d>1$ is more delicate and is related to the proof obtained in \cite{HHS} on the absence of the bound states at the threshold in the triplet $S$-sector for Schr\"odinger operators defined on a class of antisymmetric functions and where some properties of fermionic wave functions were considered. We show in Section \ref{HardyIn} that for a fixed $d$ we have $H_A(dN) \sim N^{2+ 2/d}$.

\smallskip
\noindent
Note that some related inequalities were obtained in \cite {HHLT}, \cite{L}, where it was proved that 
$$
\int_{\Bbb R^N} |\nabla u|^2\, dx \ge \frac{N^2}{4} \, \int_{\Bbb R^N} \frac{|u|^2}{|x|^2} \, dx
$$
and where $u(x) = -u(-x)\in \mathcal H^1(\Bbb R^N)$, $N\ge 2$.
If $d=1$ and $N=2$ then the constant in this inequality coincides with $C_A(2)=1$.

\smallskip
\noindent
Note that in \cite{Naz} the author has obtained Hardy's and Sobolev's  inequalities in cones. In this paper the constants in Hardy inequalities depend on the lowest Dirichlet eigenvalue of the Laplace-Beltrami operator defined on the intersection of $\Bbb S^{N-1}$ and the cone in $\Bbb R^N$. In our case such eigenvalues can be computed explicitly due to the properties of antisymmetric functions.

\smallskip
\noindent
The proof of the Sobolev  inequality \eqref{Sobolev_A} is related to a split of the space, where the antisymmetric function have the same absolute values. After that in each part one can use isoperimetric type inequalities related rearrangements.

\smallskip
\noindent
In Section \ref{degree} we study some properties of Vandermonde determinants. In Sections \ref{HardyIn} and \ref{SobolevIn} we prove Theorem \ref{Hardy_A} and Theorem \ref{Sobolev_A} respectively.
Finally we apply our results to the study of spectral properties of Schr\"odinger operators in Section \ref{Appl} and in Section \ref{Numerics} we present a table with some eigenvalues of the Laplace-Beltrami operator on antisymmetric functions defined on $\Bbb S^{dN - 1}$.


\section {\label{degree} On Vandermonde type determinants.}
\setcounter{equation}{0}

\subsection{Some preliminary results}
 Let us consider the unitary monomials of $d$ variables lexicographically, i.e. for $t\in \mathbb{R}^d$ denote $\varphi_1^{(d)}(t) = 1,$ \  $\varphi_2^{(d)}(t) = t_1$ and $\varphi_{d+1}^{(d)}(t) = t_d, \  \varphi_{d+2}^{(d)}(t) = t_1^2$  and so on. 
 We now consider the determinant
\begin{equation}\label{det}
\psi_N^{(d)}(x_1,\dots,x_N) = 
\begin{vmatrix}
  \varphi_1^{(d)}(x_1) & \cdots & \varphi_1^{(d)}(x_N)\\
  \vdots & \ddots & \vdots \\
  \varphi_N^{(d)}(x_1) & \cdots & \varphi_N^{(d)}(x_N)
\end{vmatrix}.
\end{equation}
Let $\mathcal{V}_d(N)$ be the degree of $\psi_N^{(d)}.$ 
The total degree of the above determinant is equal to the sum of degrees in every row.

\begin{proposition}\label{ppsi1}
For $N\ge 2$ the polynomial $\psi_N^{(d)}(x_1,\dots,x_N)$ is an antisymmetric, homogeneous and harmonic  function.
\end{proposition}
\begin{proof} The first two statements are obvious. In order to proof the harmonicity we use induction. Let $d$ be fixed. The base of induction is obvious, because the degree of $\psi_2^{(d)}(x_1, x_2)$ equals $1$. In the general case,
\begin{multline*}
\Delta \psi_N^{(d)}(x_1,\dots,x_N) \\
= (-1) ^N \Delta \Bigl(\varphi_N^{(d)}(x_1)\psi_{N-1}^{(d)}(x_2,\dots,  x_{N}) 
- \dots \\
\dots (-1)^N \varphi_N^{(d)}(x_N)\psi_{N-1}^{(d)}(x_1,\dots,x_{N-1})\Bigr) \\
=
\begin{vmatrix}
  \varphi_1^{(d)}(x_1) & \cdots & \varphi_1^{(d)}(x_N)\\
  \vdots & \ddots & \vdots \\
  \Delta_1\varphi_N^{(d)}(x_1) & \cdots & \Delta_N\varphi_N^{(d)}(x_N)
\end{vmatrix}.
\end{multline*}
The last equality is the consequence of the inductional assumption. It suffices to note that $\Delta_i\varphi_N^{(d)}(x_i)$ is a linear combination of $\varphi_1^{(d)}(x_i), \dots, \varphi_{N-1}^{(d)}(x_i)$ for all $1\le i \le N,$ because it has a smaller degree. The proof is complete.
\end{proof}

\begin{proposition} \label{ppsi3}
Let $P(x_1, \dots, x_N)$ be an antisymmetric homogeneus polynomial. Then $\deg P \ge \deg \psi_N^{(d)}.$
\end{proposition}
\begin{proof}
Consider $P$ as the sum of monomials
$$
P(x_1,\dots,x_N) = \sum_i a_i Q_i(x_1,\dots,x_N) = \sum_i a_i\prod_{j=1}^N q_{ij}(x_j),
$$
where $a_i\in \mathbb{C}\setminus\{0\}$ and $q_{ij}(t) = q_{ij}(t_1,\dots,t_d) = t_1^{\alpha_{ij1}}\dots t_d^{\alpha_{ijd}}$ for some $\alpha_{ijk} \in \mathbb{N}_0.$ Since $P$ is homogeneus, every $Q_i$ has the same degree as $P.$ Also for all $i$ the monomials $q_{ij}$ are pairwise distinct. Actually, if $q_{ij_1} = q_{ij_2},$ then 
$$
Q_i(x_1,\dots,x_{j_1},\dots,x_{j_2},\dots,x_n) = Q_i(x_1,\dots,x_{j_2},\dots,x_{j_1},\dots,x_n),
$$
and
$$
P(x_1,\dots,x_{j_1},\dots,x_{j_2},\dots,x_n) = - P(x_1,\dots,x_{j_2},\dots,x_{j_1},\dots,x_n).
$$
Due to equality of these polynomials we obtain the equality of respective coefficients and conclude that $a_i = 0.$  Proposition \ref{ppsi3} is proved.
\end{proof}

\noindent
Propositions \ref{ppsi1}-\ref{ppsi3} demonstrate that $\psi_N^{(d)}$ defined above is the minimal antisymmetric harmonic homogeneous polynomial.

\subsection{Exact expression of the determinant degree}

For further considerations we need to prove an auxiliary combinatorial fact.
\begin{lemma}\label{comb}
For all $n \in \mathbb{N}=\{n\}_{n=1}^\infty$ and $m \in \mathbb{N}_0 = \{0\}\cup  \mathbb{N}$ the following equality holds true
\begin{equation}\label{combin}
\sum_{k=0}^m {n + k \choose k} = \sum_{k=0}^m {n + k \choose n} = {n + m + 1 \choose n + 1}.
\end{equation}
\end{lemma}
\begin{proof}
Let us prove this proposition by using induction. For $m = 0$ the left hand side in \eqref{combin} equals ${n \choose 0} = 1$ and the right hand side becomes ${n + 1 \choose n + 1} = 1$. 
This proves the induction base. Under the assumption that \eqref{combin} is true for $m = m_0,$ it is easy to see that
$$
\sum_{k=0}^{m_0 + 1} {n + k \choose k} =  {n + m_0 + 1 \choose n + 1} + {n + m_0 + 1 \choose n} = {n + m_0 + 2 \choose n + 1}.
$$
This proves the lemma.
\end{proof}

\noindent
Let $K^{(d)}_p$ be the number of different unitary monomials of degree $p$ of $d$ variables. Note that
$$
K^{(d)}_p = \#\{\alpha = (\alpha_1, \dots, \alpha_d): \alpha_1 + \ldots + \alpha_d = p, \ \alpha_j \ge 0\}.
$$
\begin{proposition}
For $p \in \mathbb{N}_0$ and $d \in \mathbb{N}$ we have the following equality 
$$
K^{(d)}_p = {p+d-1 \choose d-1}.
$$
\end{proposition}
\begin{proof}
As previously we use induction but now with respect to $d$. It is obvious that for $d = 1$ and any $p$ we have 
$$
K^{(1)}_p = \#\{\alpha = (\alpha_1): \alpha_1 = p, \ \alpha_1 \ge 0\} = 1 = {p \choose 0}.
$$
To prove the induction step one can check that
$$
K^{(d)}_p = \sum_{\alpha_d = 0}^{p}K^{(d - 1)}_{p - \alpha_d}.
$$
Under the induction assumption we have 
$$
K^{(d)}_p = \sum_{\alpha_d = 0}^{p}{p - \alpha_d + d - 2 \choose d - 2} \underset{\beta_d = p - \alpha_d}{=} \sum_{\beta_d = 0}^{p}{\beta_d + d - 2 \choose d - 2}.
$$
According to Lemma \ref{comb}, we conclude that
$$
K^{(d)}_p = {p + d - 1 \choose d - 1},
$$
that ends the proof.
\end{proof}

For some special cases $N$ the matrix contains all monomials with degree less or equal $m$. For a fixed $m\ge0$ such $N$ we denote by $N^{(d)}_m$. Obviously,
$$
N^{(d)}_m = \sum_{p=0}^{m}K^{(d)}_p = \sum_{p=0}^{m}{p+d-1 \choose d-1} = {m+d \choose d} = \frac{(m+1)\ldots(m+d)}{d!}.
$$
The last equality is true due to Lemma \ref{comb}. Note that it implies that
$$
N^{(d)}_{m+1} = \frac{m+d+1}{m+1}N^{(d)}_m \quad {\rm and} \quad
N^{(d)}_{m+1} - N^{(d)}_m = \frac{d}{m+1}N^{(d)}_m.
$$
We now calculate the degree of $\psi_N^{(d)}.$ For the special cases 
\begin{multline*}
\mathcal{V}_d(N^{(d)}_m) = \sum_{p=0}^{m}pK^{(d)}_p = \sum_{p=1}^{m}p{p+d-1 \choose d-1} = \sum_{p=1}^{m}d{p+d-1 \choose d} \\
=d\sum_{p=0}^{m-1}{p+d \choose d} = d{m+d \choose d + 1} = \frac{d\,m}{d+1}{m+d \choose d} = \frac{d\,m}{d+1}N^{(d)}_m.
\end{multline*}
For the intermediate cases with $N^{(d)}_m < N < N^{(d)}_{m+1}$
\begin{multline*}
\mathcal{V}_d(N) = \mathcal{V}_d(N^{(d)}_m) + (N-N^{(d)}_m)(m+1) \\
= N^{(d)}_m\Bigl(\frac{dm}{d+1} - m - 1\Bigr) + N(m+1)\\
= N(m+1) - \frac{m+d+1}{d+1}N^{(d)}_m.
\end{multline*}
Note that  for $N = N^{(d)}_m$ and $N = N^{(d)}_{m+1}$ the value $\mathcal{V}_d(N)$
coincides with $\mathcal{V}_d(N^{(d)}_m)$ and $\mathcal{V}_d(N^{(d)}_{m+1})$ respectively. 
In fact, 
\begin{multline*}
N^{(d)}_{m+1}(m+1) - \frac{m+d+1}{d+1}N^{(d)}_m = N^{(d)}_{m+1}(m+1) - \frac{m+1}{d+1}N^{(d)}_{m+1} \\
=\frac{(m+1)d}{d+1}N^{(d)}_{m+1} = \mathcal{V}_d(N^{(d)}_{m+1}).
\end{multline*}
%

\subsection{Estimates of $\mathcal{V}_d(N)$}
In order to estimate $\mathcal{V}_d(N)$ we need some auxiliary facts.
Denote 
$$
A^{(d)}(m) = \frac {(m+1) + (m+2) + \ldots + (m+d)}{d},
$$  
and
$$
G^{(d)}(m) = \sqrt[d]{(m+1)\dots(m+d)}
$$
for some $m, d\in\mathbb{N}$. These are expressions for arithmetic and geometric means. It is well known that $A^{(d)}(m) \ge G^{(d)}(m)$.

\begin{lemma}
For all $m, d \in \mathbb{N}$ we have
$$
G^{(d)}(m) \ge \sqrt{(m+1)(m+d)}.
$$ 
\end{lemma}
\begin{proof}
Let $t=m + \frac{d+1}2.$ If $d$ is odd, then we have that for $k,$ s.t. $d = 2k + 1$
\begin{multline*}
\Bigl(G^{(d)}(m)\Bigr)^d = (t-k)(t - (k-1))\dots (t + k) \\
= t(t^2 - 1)(t^2 - 4)\dots (t^2 - k^2) \ge t(t^2 - k^2)^k.
\end{multline*}
According to the remark above, $G^{(d)}(m) \le t.$ Therefore
$$
\Bigl(G^{(d)}(m)\Bigr)^{2k} \ge (t^2 - k^2)^k
\quad {\rm and } \quad
G^{(d)}(m) \ge \sqrt{(m+1)(m+d)}.
$$
If $d$ is even, then we assume that $d = 2k$ and hence
\begin{multline*}
\left(G^{(d)}(m)\right)^d = \left(t-k + \frac 1{2}\right)\left(t - \left(k-\frac{3}2\right)\right)\dots \left(t + k - \frac 1{2}\right) \\
= \left(t^2 - \frac{1}4\right)\left(t^2 - \frac {9}4\right)\dots 
\left(t^2 - \left(k - \frac 1{2}\right)^2\right) \ge \left(t^2 - \left(k - \frac 1{2}\right)^2\right)^k.
\end{multline*}
\end{proof}

\begin{lemma}
For all $d\in\mathbb{N}$ there is a constant $C_d \in \mathbb{R}$ such that
$$
0 \le A^{(d)}(m) - G^{(d)}(m) < \frac{C_d}{m}.
$$
\end{lemma}
\begin{proof}
Due to the inequality on arithmetic and geometric means $A^{(d)}(m) - G^{(d)}(m) \ge 0.$ However, according to the previous lemma, 
\begin{multline*}
A^{(d)}(m) - G^{(d)}(m) \le m + \frac{d+1}2 - \sqrt{(m+1)(m+d)}
\\
=\frac{(m + \frac{d+1}2)^2 - (m+1)(m+d)}{m + \frac{d+1}2 + \sqrt{(m+1)(m+d)}} \\
= \frac{(d-1)^2}{4m + 2(d+1) + 4\sqrt{(m+1)(m+d)}}.
\end{multline*}
The last expression implies that it suffices to let $C_d$ be equal $\frac{(d-1)^2}8$.
\end{proof}

\begin{theorem} \label{deg}
Let $d \in \mathbb{N}$ and $\mathcal{V}_d(N) = \deg \psi_N^{(d)}.$ Then
$$
\mathcal{V}_d(N) = \frac{d}{d+1}\sqrt[d]{d!}\, N^{1 + \frac1{d}} - \frac{d}2 N 
+ O\left(N^{1 - \frac 1{d}}\right),
$$
as $N \rightarrow \infty.$ 
\end{theorem}
\begin{proof}
Denote $\xi_d(N) = \frac{d}{d+1}\sqrt[d]{d!}\,N^{1 + \frac1{d}} - \frac{d}2 N.$ Considering it as a function of the continuous variable $N$  we obviously have 
$\xi_d^{\prime}(N) = \sqrt[d]{d!N} - \frac{d}2$ and since $\xi_d^{\prime}$ is increasing we conclude that $\xi_d$ is convex.

\noindent
Firstly, consider the difference between $\mathcal{V}_d(N)$ and $\xi_d(N)$ at special points.
\begin{multline*}
\mathcal{V}_d(N_m^{(d)}) - \xi_d(N_m^{(d)}) = \frac{dm}{d+1}N_m^{(d)} - \frac{d}{d+1}\sqrt[d]{d!} \left(N_m^{(d)}\right)^{1 + \frac1{d}} + \frac{d}2 N_m^{(d)}  \\
= N_m^{(d)}\left(\frac{dm}{d+1} - \frac{d}{d+1}\sqrt[d]{d!N_m^{(d)}} + \frac{d}2 \right) \\
=\frac{d}{d+1}N_m^{(d)}\left(m + \frac{d+1}2 - \sqrt[d]{(m+1)\dots(m+d)}\right) \\
=\frac{d}{d+1}N_m^{(d)}\left(A^{(d)}(m) - G^{(d)}(m)\right).
\end{multline*}
Consequently, the difference $\mathcal{V}_d(N_m^{(d)}) - \xi_d(N_m^{(d)}) = o(N_m^{(d)})$ and positive. Due to the convexity of $\xi_d(N)$ and the linearity of $\mathcal{V}_d(N)$ on the segment 
$[N_m^{(d)}, N_{m+1}^{(d)}],$ the value $\xi_d(N)$ does not exceed $\mathcal{V}_d(N)$ for all $N$.

\noindent
Secondly, for all $N \ge N_m^{(d)}$ we find
$$
\xi_d(N) \ge \xi_d \left(N_m^{(d)}\right) + \left(N - N_m^{(d)}\right)\xi_d^{\prime}\left(N_m^{(d)}\right).
$$
Hence, for $N \in \left[N_m^{(d)}, N_{m+1}^{(d)}\right)$
\begin{multline*}
0 \le \mathcal{V}_d(N) - \xi_d(N) \\
\le
\mathcal{V}_d(N_m^{(d)}) + (N - N_m^{(d)})(m+1) - \xi_d(N_m^{(d)}) - \left(N - N_m^{(d)}\right)\xi_d^{\prime}(N_m^{(d)}) \\ 
= \frac{d}{d+1}N_m^{(d)}\left(A^{(d)}(m) - G^{(d)}(m)\right) + \left(N - N_m^{(d)}\right)\left(m+1 - \sqrt[d]{d!N_m^{(d)}} + \frac{d}2\right).
\end{multline*}
Finally, we note that $1 \ge \frac{N_m^{(d)}}{N} \ge \frac{N_m^{(d)}}{N_{m+1}^{(d)}} = \frac{m+1}{m+d+1}$ and consequently,
\begin{multline*}
0 \le \frac{\mathcal{V}_d(N) - \xi_d(N)}N \\
\le\frac{d}{d+1}\frac{N_m^{(d)}}N\left(A^{(d)}(m) - G^{(d)}(m)\right) + (1 - \frac{N_m^{(d)}}N)\left(A^{(d)}(m) - G^{(d)}(m) + \frac 1{2}\right)\\
\le \frac{d}{d+1}\left(A^{(d)}(m) - G^{(d)}(m)\right) + \frac{d}{m+d+1}\left(A^{(d)}(m) - G^{(d)}(m) + \frac 1{2}\right) \\
\le \frac{d}{d+1}\frac{C_d}{m} + \frac{d}{m+d+1}\left(\frac{C_d}{m} + \frac 1{2}\right) \le \frac{D_d}{m+d+1},
\end{multline*}
where $D_d = \frac{d(2d+3)}{d+1}\, C_d + \frac{d}2$. Since 
 $$
 (m + d + 1)^d > d!N_{m+1}^{(d)} > d!N,
 $$ 
 we have
 \begin{multline}\label{N-xi}
 0 \le \mathcal{V}_d(N) - \xi_d(N) \le \frac{ND_d}{\sqrt[d]{d!N}} \\
 = \frac{1}{\sqrt[d]{d!}}\Bigl(\frac{d}2 + \frac{d(d-1)^2(2d+3)}{8(d+1)}\Bigr)N^{1 - \frac 1{d}}.
 \end{multline}
\end{proof}

\begin{remark}
For $d=1$ the estimate \eqref{N-xi} is exact  $\mathcal{V}_1(N) = \xi_1(N) = \frac {N^2 - N}2$. It follows from the equality of $A^{(1)}(m)$ and $G^{(1)}(m)$. 
\end{remark}
\begin{remark}
For $d=2$ and $d=3$ we have
$$
\mathcal{V}_2(N) = \frac{2\sqrt{2}}{3}N^{\frac3{2}} - N + O(N^{\frac 1{2}}),
$$
and
$$
\mathcal{V}_3(N) = \frac{3\sqrt[3]{6}}{4}N^{\frac4{3}} - \frac{3}2 N + O(N^{\frac 2{3}}).
$$
\end{remark}

\begin{remark} \label{estmu}
Due to the proof of Theorem \ref{deg}  we conclude 
$$
\mathcal{V}_d(N) \ge \frac{d}{d+1}\sqrt[d]{d!}\, N^{1 + \frac1{d}} - \frac{d}2 N.
$$
\end{remark}

\section{Hardy inequality}\label{HardyIn}

\subsection{Laplace-Beltrami operator on $\mathbb{S}^{M-1}$.} For $M\ge2$ the Laplacian can be written in polar coordinates $(r, \theta),$ where $r = |x|$ and $\theta$ is an angular component of $x\in\mathbb{R}^M,$ in the following way
$$
-\Delta = - \frac{\partial^2}{\partial r^2} - \frac{M-1}{r}\frac{\partial}{\partial r} - \frac 1{r^2}\Delta_{\theta}.
$$
The operator $\Delta_{\theta}$ is the Laplace-Beltrami operator on $\mathbb{S}^{M-1}.$ It is well known, that the harmonic homogeneus polynomials are connected with the spherical harmonics, which are the eigenfunctions of the Laplace-Beltrami operator $-\Delta_{\theta}.$ To be more precise, let $\psi$ be the harmonic homogeneus polynomial of the degree $P$ and $\psi_{\theta} = \frac{\psi}{r^M}$. Then
$$
-\Delta_{\theta}\psi_{\theta} = P(P + M - 2)\psi_{\theta}.
$$

\begin{proposition} \label{estlambda}
Let $dN \ge 3$. The Laplace-Beltrami operator $-\Delta_{\theta}$ defined on antisymmetric functions from $L^2(\mathbb{S}^{dN - 1})$ satisfies the inequality 
$$
-\Delta_{\theta} \ge \lambda_d(N)
$$
in the quadratic form sence, where $\lambda_d(N) = \mathcal{V}_d(N)(\mathcal{V}_d(N) + Nd - 2).$
\end{proposition}
\begin{proof}
Let $\mathcal{B}$ be the orthonormal system of spherical harmonic functions and let $\mathcal{B}_A \subset \mathcal{B}$ be the orthonormal subset of the set $\mathcal{B}$ that are restrictions of antisymmetric homogeneous harmonic polynomials on $\Bbb S^{dN-1}$. For any $u \in \mathcal H^1_A(\mathbb{R}^{dN})$ we have
$$
u(r, \theta) = \sum_{k:\psi_{\theta, k} \in \mathcal{B}_A} u_k(r)\psi_{\theta, k}(\theta).
$$
According to Proposition \ref{ppsi3}, for all $k$ such that $ \psi_{\theta, k} \in \mathcal{B}_A$  
$$
\deg \psi_{k} \ge \mathcal{V}_d(N) = \frac{d}{d+1}\sqrt[d]{d!}N^{1 + \frac1{d}} - \frac{d}2 N + O \left(N^{1 - \frac 1{d}}\right).
$$ 
Hence,
\begin{multline*}
- \int_{\mathbb{S}^{dN - 1}} \Delta_{\theta}u(r, \theta) \cdot u(r, \theta) d\theta \\
= \sum_{k:\psi_{\theta, k} \in \mathcal{B}_A} |u_k(r)|^2\int_{\mathbb{S}^{dN - 1}} (-\Delta_{\theta}\psi_{\theta, k}(\theta)) \cdot \psi_{\theta, k}(\theta) d\theta  \\
\ge \sum_{k:\psi_{\theta, k} \in \mathcal{B}_A} |u_k(r)|^2 \mathcal{V}_d(N)(\mathcal{V}_d(N) + dN - 2) \\
= \lambda_d(N)\int_{\mathbb{S}^{dN - 1}} |u(r, \theta)|^2 d\theta.
\end{multline*}
\end{proof}

\begin{remark}
According to Theorem \ref{deg},
\begin{multline*}
\lambda_d(N) = \left(\frac{d}{d+1}\sqrt[d]{d!}N^{1 + \frac1{d}} - \frac{d}2 N + O(N^{1 - \frac 1{d}})\right)\\
\times\left(\frac{d}{d+1}\sqrt[d]{d!}N^{1 + \frac1{d}} - \frac{d}2 N + O(N^{1 - \frac 1{d}}) + dN - 2\right) \\ 
=\frac{d^2}{(d+1)^2}\sqrt[d]{d!^2}N^{2 + \frac2{d}} + O(N^2). 
\end{multline*}
\end{remark}
\begin{remark}
According to the Remark \ref{estmu},
\begin{multline*}
\lambda_d(N) \ge \left(\frac{d}{d+1}\sqrt[d]{d!}N^{1 + \frac1{d}} - \frac{d}2 N \right) \\
\times \left(\frac{d}{d+1}\sqrt[d]{d!}N^{1 + \frac1{d}} - \frac{d}2 N + dN - 2\right) \\
=\frac{d^2}{(d+1)^2}\sqrt[d]{d!^2}N^{2 + \frac2{d}} -\frac{d^2}4N^2 - \frac{2d}{d+1}\sqrt[d]{d!}N^{1 + \frac1{d}} + dN.
\end{multline*}
\end{remark}

\subsection{Proof of Theorem \ref{Hardy_A}. Hardy inequality for antisymmetric functions on $\mathbb{R}^{dN}$.}

Passing to polar coordinates $(r, \theta)$, we have
\begin{equation*}
\int_{\mathbb{R}^{dN}} |\nabla u (x)|^2 dx 
= \int_0^{\infty}\int_{\mathbb{S}^{dN - 1}}\left( \left|\frac{\partial u}{\partial r}\right|^2 
+ \frac{1}{r^2}\, |\nabla_\theta u|^2\right) \, r^{dN - 1} d\theta dr.
\end{equation*}
Due to the classical Hardy inequality
\begin{equation}\label{HardyClass}
\int_0^{\infty} \Bigl|\frac{\partial u}{\partial r} \Bigr|^2r^{dN-1}dr \ge \frac {(dN-2)^2}4 \int_0^{\infty} \frac{|u|^2}{r^2} r^{dN-1}dr.
\end{equation} 
Besides, Proposition \ref{estlambda} implies
\begin{equation}\label{LBhardy}
\int_{\mathbb{S}^{dN - 1}} |\nabla_\theta u|^2  \  d\theta \ge \lambda_d(N)\int_{\mathbb{S}^{dN - 1}} |u|^2 d\theta.
\end{equation}
Finally, we conclude that
\begin{multline*}
\int_{\mathbb{R}^{dN}} |\nabla u (x)|^2 dx \\
\ge \int_{\mathbb{S}^{dN - 1}}\int_0^{\infty}\left( \frac {(dN-2)^2}4\frac{|u|^2}{r^2}  + \frac {\lambda_d(N)}{r^2}|u|^2\right) r^{dN - 1} dr d\theta \\
= \left( \frac{(dN-2)^2}4+ \lambda_d(N)\right)\int_{\mathbb{R}^{dN}} \frac{|u|^2}{|x|^2} dx,
\end{multline*}
where $\lambda_d(N) = \mathcal{V}_d(N)(\mathcal{V}_d(N) + Nd - 2)$.

\begin{remark}
Using properties of $\mathcal{V}_d(N)$ we find
$$
H_A(dN) = \frac{d^2}{(d+1)^2}\sqrt[d]{d!^2}N^{2 + \frac2{d}} + O(N^2)
$$
and
\begin{multline*}
H_A(dN) \\
\ge \frac{d^2}{(d+1)^2}\sqrt[d]{d!^2}N^{2 + \frac2{d}} -\frac{d^2}4N^2 - \frac{2d}{d+1}\sqrt[d]{d!}N^{1 + \frac1{d}} + dN + \frac{(dN-2)^2}4 \\
= \frac{d^2}{(d+1)^2}\sqrt[d]{d!^2}N^{2 + \frac2{d}} - \frac{2d}{d+1}\sqrt[d]{d!}N^{1 + \frac1{d}} + 1 \\
= \left(\frac{d}{d+1}\sqrt[d]{d!}N^{1 + \frac1{d}} - 1\right)^2.
\end{multline*}
\end{remark}


\medskip
\noindent
The constant $H_A(dN)$ in Theorem \ref{Hardy_A} is sharp. Indeed,

\begin{proposition}
For $d\ge 1$ and $N \ge 2$
$$
H_A(dN) = \underset{\underset{u \not = 0}{u \in \mathcal H^1_A(\mathbb{R}^{dN})}}{\inf}
\frac{\int_{\mathbb{R}^{dN}} |\nabla u (x)|^2 dx}{\int_{\mathbb{R}^{dN}} \frac {|u(x)|^2}{|x|^2}dx}.
$$
\end{proposition}
\begin{proof} The constant $H_A(dN)$ has two terms $(dN-2)^2/4$ and $\lambda_d$. the first one is the classical Hardy constant in \eqref{HardyClass} that is sharp but not achieved. Since the inequality \eqref{LBhardy} is also sharp due to Proposition \ref{estlambda} we conclude the proof.
\end{proof}


\section{Sobolev inequality}\label{SobolevIn}

\noindent
We now consider the  Sobolev inequality on antisymmetric functions. It is well-known that for any $u\in \mathcal H^1(\mathbb{R}^{n}),$ $n \ge 3$
\[
\int_{\mathbb{R}^{n}}|\nabla u(x)|^2 dx \ge S(n) \left(\int_{\mathbb{R}^{n}}|u(x)|^{\frac{2n}{n-2}}dx\right)^{\frac{n-2}{n}},
\]
where 
$$
S(n) = \pi n(n-2)\left(\frac{\Gamma(\frac{n}2)}{\Gamma(n)}\right)^{\frac 2{n}}.
$$
The same inequality holds for any $u\in \mathcal H_0^1(\Omega),$ where $\Omega \subset \mathbb{R}^n$. Before proving our Theorem \ref{Sobolev_A} we need to study some properties of symmetric group acting on $\Bbb R^{dN}$.
\subsection{The symmetric group acting on $\mathbb{R}^{dN}$.}

Let us consider the action of symmetric group $\mathcal S_N$ on the space $\mathbb{R}^{dN}$. For an arbitrary $x = (x_1, \ldots, x_N)\in \mathbb{R}^{dN}$ and $\sigma \in \mathcal S_N$ denote by $\sigma x = (x_{\sigma(1)}, \ldots, x_{\sigma(N)})$ the permutation of elements $x$.

\begin{lemma} \label{conlem}
Let $u$ be an antisymmetric function on $\mathbb{R}^{dN}$. For an arbitrary $x \in \mathbb{R}^{dN}$ and 
$\sigma \in \mathcal S_N, \sigma \not = id$ there is no continuous path $\Gamma: [0, 1] \mapsto \mathbb{R}^{dN}$ such that $\Gamma(0) = x$, $\Gamma(1) = \sigma x$ and $|u(\Gamma(t))| > 0$ for all $t \in [0, 1].$
\end{lemma}
\begin{proof}
 We will show this by contradiction. Let $\Gamma$ be an appropriate path and $A = \{x \in \mathbb{R}^{dN} :  x_i = x_j$ for some $ 1 \le i, j \le N\}.$ Since $u(x) = 0$ for all $x \in A,$ $\Gamma([0, 1]) \cap A = \emptyset.$ Then let us consider the action $\mathcal S_N$ on $\mathbb{R}^{dN} \setminus A.$ Because of lack of fixed points on $R = \mathbb{R}^{dN} \setminus A$, this acting generates the covering 
 $p: R \mapsto  {R}/{\mathcal S_N}.$ According to the path lifting property, there are no any paths $\tilde{\Gamma}$ different from $\Gamma$ such that $p(\Gamma) = p(\tilde{\Gamma}).$

\medskip
\noindent
 Now we consider sets $E^k_{\tau} = \{x \in \mathbb{R}^{dN} :  x_{\tau(1)k} \le \ldots \le x_{\tau(N)k} \},$ where $1 \le k \le d$ and $\tau \in \mathcal S_N.$ Without loss of generality we can assume that $x \in E^k_{id}$ for all $k.$ For every $k$ we can construct the projection mapping $\Sigma^k$ on the fundamental domain $E_{id}^k$ which maps $x \in E^k_{\tau}$ to $\tau^{-1}x.$ Then if for some $k$ the point $\sigma x$ does not belong to $E_{id}^k$, then the path $\Sigma^k(\Gamma)$ differs from $\Gamma$ and $p(\Sigma^k(\Gamma)) = p(\Gamma).$ Consequently, $x$ and $\sigma x$ belong to $E_{id}^k$ for all $k$ at the same time. 

\medskip
\noindent
 It remains to note that the inequalities 
$$
x_{1k} \le \ldots \le x_{Nk},
$$
$$
x_{\sigma(1)k} \le \ldots \le x_{\sigma(N)k}
$$
imply the equality of respective parts $x_{ik} = x_{\sigma(i)k}$ for all $1 \le i \le N$ and $1 \le k \le d.$ Since $\sigma \not = id,$ there is $1 \le i \le N$ such that $\sigma(i) \not = i$ and it contradicts the assumption that $x \not \in A.$
\end{proof}

\noindent
Analogically, let us denote $\sigma E = \{\sigma x : x \in E\}$ for arbitrary $E\subset \mathbb{R}^{dN}$ and  $\sigma \in \mathcal S_N.$ Let $C_{0,A}^{\infty}(\mathbb{R}^{dN})$ be the class of $C_0^\infty$ antisymmetric functions on $\mathbb{R}^{dN}$. 

\begin{theorem} \label{conth}
Let $u \in C_{0,A}^{\infty}(\mathbb{R}^{dN})$. Then there is a set $E \subset \mathbb{R}^{dN}$ such that sets $\{\sigma E\}_{\sigma \in \mathcal S_N}$ are disjoint,  $u = 0$   for all $x \in \delta E$ and $\mu ({\rm supp}\,  u \ \setminus \ \bigcup_{\sigma \in \mathcal S_N} \sigma E) = 0$.
\end{theorem}
\begin{proof}
Let us divide ${\rm supp}\,  u$ into the union of connected components $\{U_{\alpha}\}$. Because of their openness, $U_{\alpha}$ are path-connected and according to Lemma \ref{conlem}, we can consider the group action $\mathcal S_N$ on $\{U_{\alpha}\}$. At the end, it remains to choose representatives from each equivalence class. Their union will be the desired set $E$. 
\end{proof}

\subsection{Proof of Theorem \ref{Sobolev_A} (Sobolev Inequality)}

It is enough to show it for an arbitrary $u$ belonging to the subclass 
$C_{0,A}^\infty(\mathbb{R}^{dN})$ of $C_0^\infty(\mathbb{R}^{dN})$-functions satisfying antisymmetry conditions. The case $u\in \mathcal H_A^1(\mathbb{R}^{dN})$ follows by the completness $C^{\infty}_{0,A}(\mathbb{R}^{dN})$ in $\mathcal H_A^1(\mathbb{R}^{dN})$. 

Let $E$ be the set from Theorem \ref{conth}.  The restriction of $u$ to the set $E$ satisfies zero boundary conditions at the boundary $\partial E$. Thus we have
$$
\int_{E}|\nabla u(x)|^2 dx 
\ge   S(dN) \left(\int_{E}|u(x)|^{\frac{2dN}{dN-2}}dx\right)^{\frac{dN-2}{dN}}.
$$
This implies 
\begin{multline*} 
\int_{\mathbb{R}^{dN}}|\nabla u(x)|^2 dx = N! \int_{E}|\nabla u(x)|^2 dx \\
\ge  N! \ S(dN) \left(\int_{E}|u(x)|^{\frac{2dN}{dN-2}}dx\right)^{\frac{dN-2}{dN}} \\
= (N!)^{\frac 2{dN}} S(dN)  \left(\int_{\mathbb{R}^{dN}}|u(x)|^{\frac{2dN}{dN-2}}dx\right)^{\frac{dN-2}{dN}}
\end{multline*} 
which proves Theorem \ref{Sobolev_A}.

\begin{remark}
 The constant $S_A(dN) = (N!)^{\frac 2{dN}} S(dN)$ is sharp and substantially larger than the classical one. It is enough to consider the minimizing sequence for the classical Sobolev inequality on $E_{id}^1$ and extend it on $\mathbb{R}^{dN}$ by antisymmetry.
\end{remark}
\begin{remark}
Due to Stirling's approximation we find
$$
S_A(dN) \sim \frac{\pi e^{1 - \frac 2{d}}}2 dN^{1 + \frac 2{d}}.
$$
\end{remark}

\section{Applications to spectral inequalities}\label{Appl}

Having two improved classical inequalities (Hardy and Sobolev) we now apply them to spectral properties of Schr\"odinger operators.

\subsection{Caffarelli-Kohn-Nirenberg type inequality.}

\begin{proposition}
Let $p = \frac{2dN}{dN - 2\nu}$ and $\gamma = 2dN\frac{\nu - 1}{dN - 2\nu}, 0 \le \nu \le 1, dN \ge 3,$    Then for any function $u \in \mathcal H^1_A(\mathbb{R}^{dN})$ we have
\begin{multline*}
\left(\int_{\mathbb{R}^{dN}}|x|^{\gamma} |u(x)|^p dx\right)^{\frac{2}{p}} \\
\le\widetilde{\mathcal K}(dN, \nu) \left(\int_{\mathbb{R}^{dN}}|\nabla u(x)|^2dx\right)^{\nu}\left(\int_{\mathbb{R}^{dN}}\frac{|u(x)|^2}{|x|^2}dx\right)^{1 - \nu},
\end{multline*}
where $\widetilde{\mathcal K}_d(N, \nu) = S_A^{-\nu}(dN).$
\end{proposition}

\begin{proof}
\begin{multline*}
\left(\int_{\mathbb{R}^{dN}}|x|^{\gamma} |u(x)|^p dx\right)^{\frac{2}{p}} = \left(\int_{\mathbb{R}^{dN}}\left(\frac{|u(x)|}{|x|}\right)^{(1 - \nu)p} |u(x)|^{p\nu} dx\right)^{\frac{2}{p}} \\
\le \left(\int_{\mathbb{R}^{dN}}|u(x)|^{\frac{2dN}{dN-2}}dx\right)^{\frac{\nu(dN-2)}{dN}}
\left(\int_{\mathbb{R}^{dN}}\left(\frac{|u(x)|}{|x|}\right)^{2}dx\right)^{1 - \nu} \\
\le S_A^{-\nu}(dN) \left(\int_{\mathbb{R}^{dN}}|\nabla u(x)|^2dx\right)^{\nu}\left(\int_{\mathbb{R}^{dN}}\frac{|u(x)|^2}{|x|^2}dx\right)^{1 - \nu}.
\end{multline*}
\end{proof}

\begin{proposition}
Let $p = \frac{2dN}{dN - 2\nu}$ and $\gamma = 2dN\frac{\nu - 1}{dN - 2\nu}, 0 \le \nu \le 1, dN \ge 3.$ Then for any antisymmetric function $u \in \mathcal H_A^1(\mathbb{R}^{dN})$ we have
$$
\int_{\mathbb{R}^{dN}}|\nabla u(x)|^2dx \ge \mathcal K_d(N, \nu)\left(\int_{\mathbb{R}^{dN}}|x|^{\gamma} |u(x)|^p dx\right)^{\frac{2}{p}},
$$
where $\mathcal K_d(N, \nu) =S_A^{\nu}(dN)\, H_A^{1 - \nu}(dN).$
\end{proposition}
\begin{proof}
\begin{multline*}
\left(\int_{\mathbb{R}^{dN}}|x|^{\gamma} |u(x)|^p dx\right)^{\frac{2}{p}} \\
\le S_A^{-\nu}(dN) \left(\int_{\mathbb{R}^{dN}}|\nabla u(x)|^2dx\right)^{\nu}\left(\int_{\mathbb{R}^{dN}}\frac{|u(x)|^2}{|x|^2}dx\right)^{1 - \nu} \\
\le
\frac 1{S_A^{\nu}(dN)H_A^{1 - \nu}(dN)}\int_{\mathbb{R}^{dN}}|\nabla u(x)|^2dx.
\end{multline*}
\end{proof}

\begin{remark} Note that as $N\to\infty$ we have 
$$
\mathcal K_d(N, \nu) \sim \Bigl(\frac{\pi e^{1 - \frac 2{d}}}2 \Bigr)^{\nu}  \frac{d^{2-\nu}}{(d+1)^{2-2\nu}}\sqrt[d]{d!^{2-2\nu}}N^{2 + \frac{2}{d} - \nu}$$
\end{remark}


\medskip
\subsection{Spectral properties of Schr\"odinger operators.}
Let us consider a Schr\"odinger operator defined on antisymmetric functions in $L_A^2(\mathbb{R}^{dN})$
$$
\mathscr{H} = -\Delta - V,
$$
where $V \ge 0.$ 

\begin{theorem}
Let $dN \ge 3$ and $0 \le \nu \le 1.$ Assume that
$$
\left(\int_{\mathbb{R}^{dN}} V^{\frac {dN}{2\nu}}|x|^{\frac{1-\nu}{\nu}dN}dx\right)^{\frac{2\nu}{dN}} \le \mathcal K_d(N, \nu).
$$
Then the operator $\mathscr{H}$ is positive and has no negative eigenvalues.
\end{theorem}
\begin{proof}
The quadratic form of $\mathscr{H}$ equals
\begin{align*}
\int_{\mathbb{R}^{dN}} &|\nabla u(x)|^2dx - \int_{\mathbb{R}^{dN}} V(x)|u(x)|^2 dx\\
&= \int_{\mathbb{R}^{dN}}|\nabla u(x)|^2dx - \int_{\mathbb{R}^{dN}} V(x)|x|^{2(1 - \nu)}|u(x)|^2|x|^{2(\nu - 1)} dx\\
&\ge \int_{\mathbb{R}^{dN}}|\nabla u(x)|^2dx - \left(\int_{\mathbb{R}^{dN}} V(x)^{\frac{dN}{2\nu}}|x|^{\frac{1-\nu}{\nu}dN}dx \right)^{\frac{2\nu}{dN}}\\
&\times \left(\int_{\mathbb{R}^{dN}}|u(x)|^{\frac{2dN}{dN-2\nu}}|x|^{2dN\frac{\nu - 1}{dN-2\nu}}dx\right)^{\frac{dN-2\nu}{dN}} dx\\
&\ge \mathcal K_d(N, \nu)\left(\int_{\mathbb{R}^{dN}}|x|^{\gamma} |u(x)|^p dx\right)^{\frac{2}{p}}\\
&- \mathcal K_d(N, \nu)\left(\int_{\mathbb{R}^{dN}}|u(x)|^{p}|x|^{\gamma}dx\right)^{\frac{2}{p}} dx = 0,
\end{align*}
where $p = \frac{2dN}{dN - 2\nu}$ and $\gamma = 2dN\frac{\nu - 1}{dN - 2\nu}.$ Consequently, $\int_{\mathbb{R}^{dN}}|\nabla u(x)|^2dx - \int_{\mathbb{R}^{dN}} V(x)|u(x)|^2 dx \ge 0$ and this completes the proof.
\end{proof}

\section{Some numerical values}\label{Numerics}

\noindent
Despite of complexity of exact values of $\mathcal{V}_d(N)$ we can give some values of minimal eigenvalues of the Laplace-Beltrami operator obtained with numerics.

\begin{center}
    \begin{tabular}{|c||*{8}{c|}}\hline
&\makebox[2.2em]{N = 2}&\makebox[2.2em]{N = 3}&\makebox[2.2em]{N = 4}&\makebox[2.2em]{N = 5}&\makebox[2.2em]{N = 6}&\makebox[2.2em]{N = 7}&\makebox[2.2em]{N = 8}&\makebox[2.2em]{N = 9}\\
\hline\hline
\makebox[2.8em]{d = 1} & 1 & 12 & 48 & 130 & 285 & 546 & 952 & 1548  \\\hline
\makebox[2.8em]{d = 2} & 3 & 12 & 40 & 84 & 144 & 253 & 392 & 561 \\\hline
\makebox[2.8em]{d = 3} & 5 & 18 & 39 & 90 & 161 & 252 & 363 & 494 \\\hline
\makebox[2.8em]{d = 4} & 7 & 24 & 51 & 88 & 168 & 272 & 400 & 552 \\\hline
\makebox[2.8em]{d = 5} & 9 & 30 & 63 & 108 & 165 & 280 & 423 & 594 \\\hline
\makebox[2.8em]{d = 6} & 11 & 36 & 75 & 128 & 195 & 276 & 432 & 620\\\hline
\end{tabular}
\end{center}

\medskip
\centerline{Minimal eigenvalues of the Laplace-Beltrami operator on $\Bbb S^{dN-1}$ }
\centerline{ on antisymmetric functions.}

\bigskip
\noindent
Also it allows us to compare $\mathcal{V}_d(N)$ and its estimate $\xi_d(N)$ that we obtained. Results are performed on following graph, where the reader can see the difference  $\mathcal{V}_d(N) - \xi_d(N)$ with growth $N$ from $2$ to $100$ for $2 \le d \le 8$:

\begin{figure}[h]
\caption{The values of $\mathcal{V}_d(N) - \xi_d(N)$, $2 \le d \le 8$}
\centering
\includegraphics[width = 13cm]{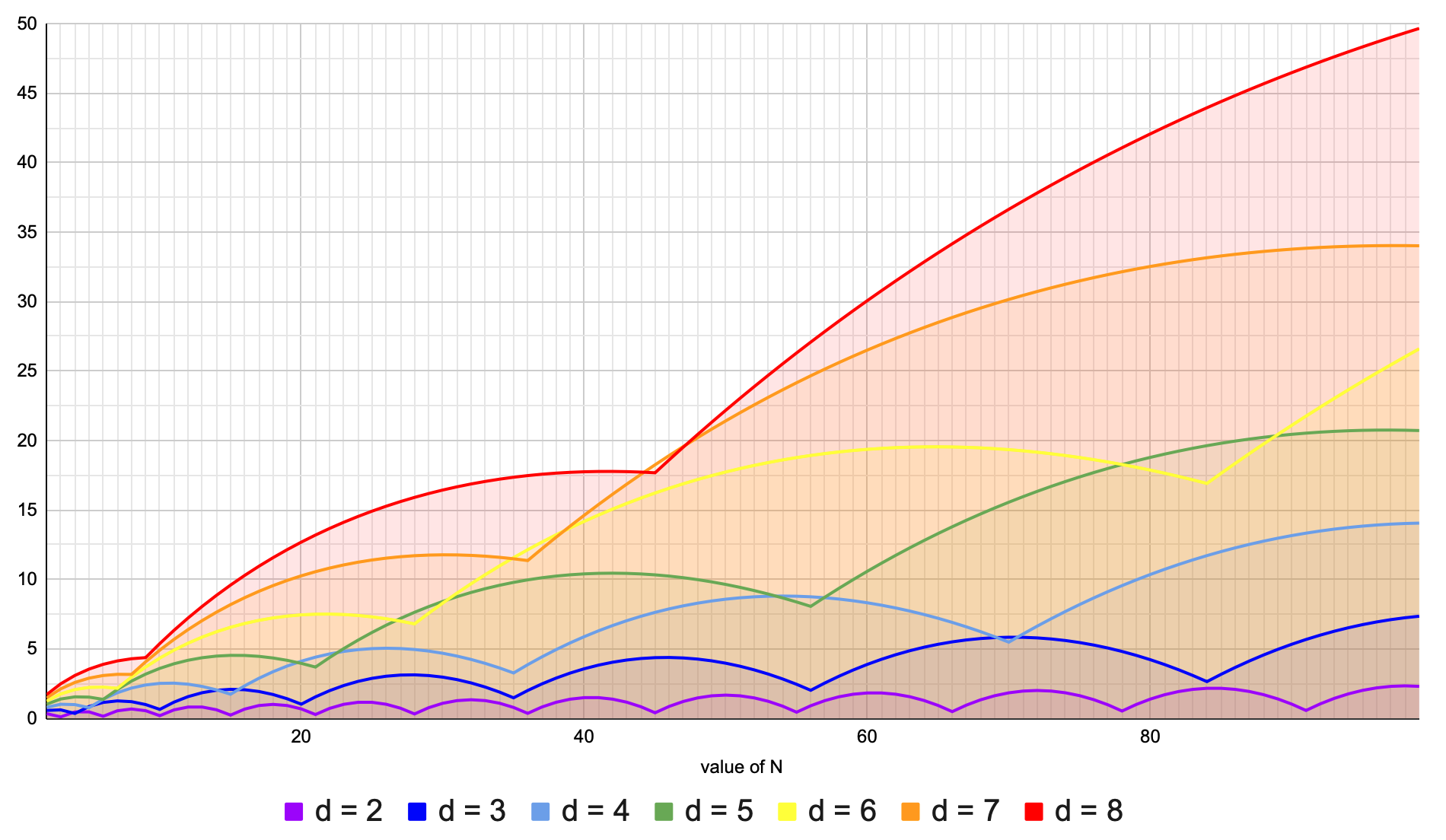}
\end{figure}

\medskip
\noindent
Due to \eqref{N-xi} the difference of $\mathcal{V}_d(N)$ and $\xi_d(N)$ equals $O(N^{1 - \frac 1{d}})$. It justifies the growth of the difference with the increasing of $d$. Also we can see the special values $N_m^{(d)}$. They correspond to cusps on graphs.

\bigskip

\noindent
{\it Acknowledgements:}  AL was supported by the Ministry of Science and Higher Education of the Russian Federation, (Agreement 075-10-2021-093, Project MTH-RND-2124).


\end{document}